\documentstyle[11pt,amssymb,diagrams]{amsart}
\newarrow{Equal}=====

\newcommand{\exseq}[3]{0 \ar #1 \ar #2 \ar #3 \ar 0}
\newcommand{\seq}[3]{{#1}_{#2}, \ldots, {#1}_{#3}}

\newcommand{\surj}{\twoheadrightarrow}
\newcommand{\inc}{\hookrightarrow}
\newcommand{\ar}{\rightarrow}
\newcommand{\bd}{\partial}
\newcommand{\x}{\times}
\newcommand{\ox}{\otimes}
\newcommand{\iso}{\cong}

\newcommand{\Gr}{\text{Gr}}
\newcommand{\Pic}{\text{Pic}\>}
\newcommand{\pn}{{\Bbb P}^n}
\newcommand{\pnr}{{\Bbb P}^n_{\Bbb R}}
\newcommand{\drk}{{D^r_k}}
\newcommand{\crk}{{C^r_k}}

\newcommand{\Hom}{\text{Hom}}

\newcommand{\reg}{\hbox{\scriptsize reg}}

\newcommand{\cF}{{\cal F}}

\newcommand{\CC}{{\Bbb C}}

\newcommand{\PP}{{\Bbb P}}

\newcommand{\RR}{{\Bbb R}}

\newcommand{\ZZ}{{\Bbb Z}}

\newcommand{\g}{\gamma}

\renewcommand{\l}{\lambda}

\newcommand{\s}{\sigma}

\renewcommand{\S}{\Sigma}

\renewcommand{\L}{\Lambda}

\theoremstyle{plain}
\begingroup
\theorembodyfont{\sl}
\newtheorem{thm}{Theorem}
\newtheorem{cor}[thm]{Corollary}
\newtheorem{lem}[thm]{Lemma}
\newtheorem{prop}[thm]{Proposition}

\endgroup

\theoremstyle{definition}
\newtheorem{defn}[thm]{Definition}

\theoremstyle{remark}
\newtheorem{rem}[thm]{Remark}
\newtheorem{ex}[thm]{Example}

\title[Degree of the divisor of solutions of a differential equation]{Degree of the
divisor of solutions of a differential equation on a projective variety}

\author{Vicente Mu\~noz and Ignacio Sols}
\thanks{Mathematics Subject Classification. Primary: 14N10. Secondary: 14D20, 14M15, 14P05.}
\date{December, 1998}

\begin{document}

\maketitle

\begin{abstract}
  Using the data schemes from~\cite{ASS} we give a rigorous definition of algebraic
  differential equations on the complex projective space $\pn$. For an algebraic subvariety
  $S \subseteq \pn$,
  we present an explicit formula for the degree of the divisor of solutions of 
  a differential equation on $S$ and give some examples of applications. 
  We extend the technique and result to the real case.
\end{abstract}

\section{Introduction}
\label{sec:1}

We deal with the problem of finding the degree of the divisor of solutions of a 
differential equation on a projective variety,
which was studied by Halphen in~\cite{H} for differential 
equations on plane curves and on hypersurfaces of $\pn$. With the notion of infinitesimal 
data in $\pn$, introduced in the plane by A. Collino~\cite{C}, used by
S. Colley and G. Kennedy in~\cite{CK1}~\cite{CK2}, and generalised to higher 
dimensions in~\cite{ASS}, we solve this problem for algebraic subvarieties of $\pn$ 
(the method of Halphen seems to generalise only to complete intersections). 
The plane data of Collino have been rediscovered by Mohamed Belghiti~\cite{Mo} and
applied to give a modern proof of the Halphen formula for plane curves, including some
explicit calculation of the involved invariants of the equation that Halphen was able
to obtain in this case.

Let $\pn$ the complex projective space of dimension $n$ (although we can work over 
any algebraically closed field $K$ of zero characteristic).
Associated to a set of coordinates $x_1,\ldots, x_k, y_1,\ldots, y_{n-k}$ on an affine
open set $U^0=\CC^n \subseteq \pn$,
we construct inductively an open set $U^r$ of the data scheme $\drk \pn$, where the
partial derivatives of the $y_j$ with respect to the $x_i$ are understood as the
canonical coordinates of $U^r$. In this context, a differential equation 
  $$
    f(x_i,y_j, {\bd y_j \over \bd x_i}, {\bd^2 y_j \over \bd x_{i_2} \bd x_{i_1}}, \cdots, 
  {\bd^r y_j \over \bd x_{i_r} \cdots \bd x_{i_1}})=0,
  $$ 
is understood as a algebraic equation on the open set $U^r$. The data satisfying the differential
equation form a divisor $\drk \pn (f)$ on $\drk \pn$ defined as the closure of the solutions to $f=0$. 
For the differential equation $f$ there are naturally well defined enumerative invariants 
$\g^f_s$, $0\leq s\leq r$.

Let $S\subseteq \pn$ be an algebraic subvariety of dimension $k$, then $S(f)$ will be the 
locus of the solutions of the differential equation $f=0$ on $S$. The main result
of the paper (theorem~\ref{thm:main}) gives a closed formula
for $\deg S(f)$ (assuming $S(f)$ is a proper subset of $S$), 
expressed as the scalar product between the $\g^f_s$ and the
cuspidal numbers $\g^s_S$ of $S$. This reduces to a simple
formula $\deg S(f)=\g_0^f\g^0_S + \g_1^f\g^1_S$
when $S$ is smooth (or normal), where $\g^0_S$ and $\g^1_S$ are
the degree and class of $S$ respectively, and $\g_0^f$ and $\g^1_f$
can be computed directly, either with the use of proposition~\ref{prop:7}
or remark~\ref{rem:util}
(when $f$ is distinguished in suitable variables, 
which is the general case) or by
using the formula $\deg S(f)=\g_0^f\g^0_S + \g_1^f\g^1_S$
applied to a smooth quadric and a smooth cubic subvarieties 
$S \subseteq \pn$, and solving the system of linear equations thus obtained.

The use of theorem~\ref{thm:main} is two-fold. On the one hand it gives the degree of $S(f)$ if
we know enough information from $f$ and $S$. This can be used to understand geometric properties 
of $S$ as these are usually described
by a suitable differential equation (e.g.\ flexes of plain curves, parabolic points).
On the other hand the explicit computations of $\deg S(f)$ in particular cases may lead
to finding the invariants $\g_s^f$ of a differential equation $f$, and the cuspidal numbers
$\g^s_S$ of a subvariety $S$.
Section~\ref{sec:4} is devoted to examples in which these two applications
of theorem~\ref{thm:main} are worked out. In particular, we find the degree of the
divisor of parabolic points of a subvariety of $\pn$.

Finally in the last section we extend the results for the real projective
space $\pnr$, by studying the behaviour of the constructions involved
via the conjugation involution. This is 
merely an introduction and not a thorough analysis of the use of data varieties in the
real case. The formula thus obtained  in theorem~\ref{thm:main2} is analogue to the 
one in the complex case, but we have to restrict ourselves to coefficients in $\ZZ/2\ZZ$.
We end up with a simple application on the number of umbilical points of a real surface in
$\PP_{\RR}^3$.

\section{Differential equations on projective space}
\label{sec:2}

To make precise sense of a differential equation in the projective space
$\pn$, we need to recall from~\cite{ASS} the definition of the smooth compact
moduli $\drk \pn$ of infinitesimal data of dimension $k$ at order $r$ in $\pn$.

Let $Z$ be a scheme smooth over $\CC$. We define schemes $\drk Z$ together with
embeddings and projections
$$
\begin{diagram}
  \drk Z &\rInto &\Gr_k T D^{r-1}_k Z \\
   \dOnto_{b_r} & \ldOnto_{\pi_r} \\
   D^{r-1}_k Z
\end{diagram}
$$
inductively, as follows. We take $D^0_k Z=Z$, $D^1_k Z=\Gr_k TZ$ and if $r\geq 2$ and
$$
\begin{diagram}
  D^{r-1}_k Z&\rInto &\Gr_k T D^{r-2}_k Z \\
   \dOnto_{b_{r-1}} & \ldOnto_{\pi_{r-1}} \\
   D^{r-2}_k Z
\end{diagram}
$$
is already defined then $\drk Z \subseteq \Gr_k TD^{r-1}_k Z$ is the Grassmannian 
$\Gr_k \cF_{r-1}$ of the (Semple) bundle $\cF_{r-1} \subseteq TD^{r-1}_k Z$ obtained as 
pull-back
\begin{equation}
\label{eqn:sols1}
\begin{diagram}
    &&&&\L_{r-1} &\rEqual& \L_{r-1} \\
     &&&& \uOnto &&  \uOnto \\
0 &\rTo&  TD^{r-1}_k Z/D^{r-2}_k Z &\rTo& TD^{r-1}_k Z &\rTo^{Tb_{r-1}}&
                                                     b_{r-1}^* TD^{r-2}_k Z &\rTo &0 \\
   &&  \dEqual &&      \uInto & \text{pull-back}&     \uInto \\
0 &\rTo & TD^{r-1}_k Z/D^{r-2}_k Z& \rTo &\cF_{r-1} &\rTo& \S_{r-1}& \rTo &0
\end{diagram}
\end{equation}
where the right hand sequence is the restriction to $D_k^{r-1} Z$ of the universal
sequence on $\Gr_k T D_k^{r-2} Z$.

On each $\drk Z =\Gr_k \cF_{r-1}$ there is a natural divisor $\crk Z$, 
namely the Schubert special cycle of $k$-planes of $\cF_{r-1}$ meeting 
$TD^{r-1}_k Z/D^{r-2}_k Z$, whose elements are called cuspidal. The
elements in the complement are called non-cuspidal.

If $Z'$ is a subvariety of $Z$ of dimension $k$ and $Z'_{\reg}$ is the open set
of smooth points, then
$$
   Z'_{\reg} \iso \drk Z'_{\reg} \subseteq \drk Z
$$
by the obvious functoriality of our construction. We then define $\drk Z' \subseteq
\drk Z$ as the closure of this subset.

From now on $Z=\pn=\PP(V)$ and $1 \leq k \leq n-1$. 
Fix a trivialisation $V \iso \CC^{n+1}$ or projective reference, with corresponding hyperplane
$C^0$ at infinity and affine part $U^0=
\pn - C^0=\CC^n =\CC^k \x \CC^{n-k}$, with affine coordinates which we name
\begin{equation}
  x_1, \ldots, x_i, \ldots, x_k,y_1, \ldots,y_j, \ldots, y_{n-k}.
\label{eqn:2}
\end{equation}
In order to make sense of a differential equation, we need to introduce formal
partial derivative symbols of the variables $y_j$ with respect to the variables 
$x_i$. For any $s \geq 1$, $1 \leq j \leq n-k$, 
and $1 \leq \seq{i}{1}{s} \leq k$ we introduce
\begin{equation} 
\label{eqn:1}
   y^j_{i_1,\ldots, i_s}={\bd^s y_j \over \bd x_{i_s} \cdots \bd x_{i_1}}.
\end{equation}
At the present stage these are to be understood merely as symbols. Later we
will identify them as coordinates of an affine open set of 
the data scheme $\drk \pn$. 
Note that we have introduced non-commuting derivative symbols 
(see remark~\ref{rem:commuting} below).

On each point $P \in U^0$ we consider the $k$-space $X(P)=P +\CC^k \subseteq TU^0$ and 
the $(n-k)$-space $Y(P)=P +\CC^{n-k} \subseteq TU^0$
providing subbundles $X$, $Y$ of $TU^0=X  \oplus Y$, 
all of them trivial. The complement $U^1$ in $D^1_k U^0=\Gr_k TU^0$ of the 
special Schubert cycle of $k$-subbundles of $TU^0$ meeting $Y$ is thus
\begin{equation} 
  U^1=\Hom (X,Y) = U^0 \x \Hom (\CC^k ,\CC^{n-k})=U^0 \x V^1,
\label{eqn:3}
\end{equation}
where $V^1=\Hom (\CC^k ,\CC^{n-k})$ 
is a cartesian power of $\CC$ with coordinates $y^j_i$ (standing for
${\bd y_j \over \bd x_i}$). 
On $U^1$, the lifted $k$-subbundle $X_{U^1} \subseteq (TU^0)_{U^1}$ provides a split of the
universal sequence 
$$
 \exseq{\S_1}{(TU^0)_{U^1}}{\L_1}
$$
on $U^1$. The universal subbundle $\S_1$ on $U^1$ becomes then isomorphic to the bundle
$X_{U^1} \iso U^1 \x \CC^k$. 

We shall denote $C^1$
the closure in $D^1_k \pn$ of the complement of $U^1$ in $D^1_k U^0$.
Also for all $r \geq 2$ we let $C^r=C^r_k \pn$ be the cuspidal divisor
of $D^r_k \pn$. Denote by $C^r_s$ for all $0\leq s\leq r$ the divisor of $\drk \pn$
counterimage by 
\begin{equation}
 b_{r,s}=b_{s+1}\circ \cdots\circ b_{r-1}\circ b_r: \drk \pn \to D^s_k \pn
\label{eqn:sols2}
\end{equation}
of the divisor $C^s$ of $D^s_k \pn$ (in particular $C^r_r=C^r$). 
Denote by $U^r$ the open subset of $\drk \pn$
$$
  U^r=\drk \pn \backslash (C^r_0 \cup \ldots \cup C^r_s\cup \ldots\cup C^r_r)
$$
consisting of data which are not in $C^r_0$ (i.e.\ which are ``finite''),
are not in $C^r_1$ (i.e.\ which are nowhere ``vertical'') and are not
in $C^r_s$ for any $s\geq 2$ (i.e.\ which are ``nonsingular'').
It is important to note that $C^r_s$ is intrinsically defined for $s\geq 2$, but
is dependent on the choice of trivialisation for $s=0,1$.
The data in $U^r$ are the 
data for which the symbols~\eqref{eqn:1}, for $s=1,\ldots, r$,
are going to acquire a precise meaning.

\begin{lem}
\label{lem:2}
  $C^r_0, C^r_1, \ldots,C^r_s, \ldots, C^r_r$ is a basis for $A^1 \drk \pn$.
\end{lem}

\begin{pf}
  We can see this by induction on $r \geq 0$. For $r=0$ it is obvious. In general, we have
  to show that the natural map 
  \begin{equation}
         \Pic D^{r-1}_k \pn \x \ZZ C^r \ar \Pic \drk \pn =\Pic \Gr_k (\cF_{r-1})
  \label{eqn:sols3}
  \end{equation} 
  is an isomorphism. We need only to remark that for
  any point $P\in  D^{r-1}_k \pn$, $(C^r)_P$ is a basis of the Picard group of the
  Grassmannian $\Gr_k((\cF_{r-1})_P)$. Surjectivity of~\eqref{eqn:sols3} now follows easily 
  from~\cite[Ex. III.12.4]{Ha} and injectivity by restricting to a fibre.

  Alternatively, the lemma follows from the explicit descriptions given below. The complement
  of the union of all $C^r_s$ is a cartesian power of $\CC$ (proposition~\ref{prop:triv}) and
  $C^r_s$ are all irreducible, therefore they generate $A^1\drk\pn$.
  On the other hand, by 
  the proof of step 2 of proposition~\ref{prop:7}, 
  they are linearly independent. So they form a basis for $A^1\drk\pn$.
\end{pf}

\begin{prop}
\label{prop:triv}
  The given trivialisation~\eqref{eqn:2} of $U^0$ induces trivialisations 
  $$U^r \iso U^{r-1} \x \Hom( (\CC^k)^{\ox r}, \CC^{n-k})$$
  of $U^r$ with cartesian powers of $\CC$. The coordinates of $U^r$ are given by~\eqref{eqn:2}
  and~\eqref{eqn:1}, $1\leq s\leq r$.
\end{prop}

\begin{pf}
  Starting from the trivialisation of $U^0$, we get an induced
  trivialisation~\eqref{eqn:3} of $U^1$. Suppose we have already trivialisations 
  $U^s=U^{s-1} \x V^s$
  with $V^s \iso \Hom((\CC^k)^{\ox s}, \CC^{n-k})$, which is isomorphic to a cartesian
  power of $\CC$, for $1 \leq s \leq r$. 
  Then $(\S_r)_{U^r} \iso (\S_{r-1})_{U^r} \iso \cdots  \iso(\S_1)_{U^r} \iso X_{U^r} 
  \iso U^r \x \CC^k$. 
  The Semple sequence on $U^r$ is split by
$$
  \cF_r \inc TU^r \surj TU^r/U^{r-1}
$$
  for all $r$, so that $\cF_r=\S_r \oplus TU^r/U^{r-1}$. Consequently,
$$
  U^{r+1} \iso \Hom ( (\S_r)_{U^r}, TU^r/U^{r-1})= U^r \x \Hom( \CC^k, V^r),
$$
  i.e. $U^{r+1} \iso U^r \x V^{r+1}$ with $V^{r+1} \iso  \Hom( (\CC^k)^{\ox (r+1)}, \CC^{n-k})$,
  as required.

  As for naming the coordinates corresponding to this chart $U^r$ of $\drk \pn$, 
  these are given by~\eqref{eqn:2} for $U^0$, $y^j_i$ for $V^1$, and in general, assuming
  $y^j_{i_1,\ldots,i_r}$ are coordinates for
  $V^r$, then the space $V^{r+1} \iso \Hom( \CC^k, V^r)$
  is described by coordinates which we denote $\bd y^j_{i_1,\ldots,i_r} \over \bd x_{i_{r+1}}$,
  i.e. $y^j_{i_1,\ldots,i_r,i_{r+1}}$.
\end{pf}

For later use, there is a zero section $0_r:U^{r-1} \ar U^{r}$ of
$b_{r}:U^{r} \ar U^{r-1}$ and composing, a zero section
$0_{s,r}:U^s \ar U^r$ of $b_{r,s}$, for all $s \leq r$.
Now we are in the position of defining differential equations on $\pn$.
\begin{defn}
\label{def:1}
  A differential equation on $\pn$ relative to a reference $x_i, y_j$ is a
  non-zero algebraic equation
  $$
    f(x_i,y_j, y^j_i, y^j_{i_1,i_2}, \ldots, y^j_{i_1,\ldots,i_r})=0,
  $$ 
  on $U^r$, for some $r \geq 1$. 
\end{defn}

  Let $U^r(f)$ be the divisor of $U^r$ defined by the differential equation $f$, and let 
  $H_f=\drk \pn (f)$ denote its closure in $\drk \pn$, i.e.\ the data
  satisfying the differential equation. In this way, any differential equation $f$
  gives a hypersurface $H_f$ in $\drk \pn$ not containing any of $C^r_2, \ldots, C^r_r$.
  Conversely, given a hypersurface $H$ in $\drk \pn$ not containing any of 
  $C^r_2, \ldots, C^r_r$, we can find a suitable reference~\eqref{eqn:2} such that
  $H$ does not contain $C^r_0$ and $C^r_1$, thus being defined as $\drk \pn (f)$ for a suitable
  polynomial $f$ on $U^r$. So we can give a more intrinsic definition

\begin{defn}
\label{def:2}
  A differential equation on $\pn$ is a hypersurface $H \subseteq \drk \pn$ not 
  containing any of $C^r_2, \ldots, C^r_r$.
\end{defn}

  This second viewpoint allows us to forget coordinates. Nonetheless, in 
  practice, we need to work with coordinates. Whenever we say that $f$ is a differential 
  equation on the projective space, we shall understand that there is a projective 
  reference implicit. There are two caveats. First, if we change the reference,
  the polynomial $f$ will change in general. Second, there might be some forbidden references
  that we cannot choose (namely, those in which $H_f$ contains either of $C^r_0$, $C^r_1$).

\begin{rem}
\label{rem:commuting}
The formal partial derivative symbols we have defined are non-commuting. This simply
means that, for instance, $y^j_{i_1,i_2}$ and $y^j_{i_2,i_1}$  ($i_1 \neq i_2$) are 
independent coordinates. We can solve this difficulty by restricting to the subset of $U^r$
given by the equations $y^j_{i_1,\ldots, i_s}=y^j_{\s(i_1),\ldots, \s(i_s)}$,
$1 \leq j\leq n-k$, $1\leq s\leq r$, $\s$ any permutation of indices, and considering the 
closure of such subset in $\drk \pn$.

More intrinsically, this subvariety of symmetric data $S^r_k Z\subseteq \drk Z$ has been 
extracted in~\cite{ASS}, for any smooth scheme $Z$. 
In our case, dealing with $\drk Z$ is simpler and will suffice. 
\end{rem}

  Associated to a differential equation we have naturally defined enumerative invariants,
  thanks to lemma~\ref{lem:2}.
\begin{defn}
\label{def:3}
 For a differential equation $f$ on $\pn$, define $\g^f_0, \ldots, \g^f_r \in\ZZ$ by 
 $$
   [\drk\pn(f)]= \g_0^f C^r_0 +\cdots +\g_r^f C^r_r 
 $$
 in $A^1 \drk\pn$.
\end{defn}

It is important to note that $\g_s^f$ do not depend on the trivialisation, as the divisor
classes of $C^r_0$ and $C^r_1$ are independent of the trivialisation as well.
The following proposition allows us to compute some of $\g^f_s$ in 
the general case.

\begin{prop}
\label{prop:7}
  Let $f$ be a differential equation on $\pn$ relative to a reference $x_i,y_j$.
\begin{itemize}
\item  Suppose $f$ is distinguished in the variable $y^1_{1,\stackrel{(r)}{\ldots},{1}}$ (i.e.\ for $d=\deg(f)$
  the monomial $(y^1_{1,\stackrel{(r)}{\ldots},{1}})^d$ appears in $f$ with nonzero coefficient), then
  the leading $\g^f_r$ is
  $$\g^f_r= \deg f(0, \ldots, y^1_{1,\stackrel{(r)}{\ldots},{1}} ,\ldots, 0).$$
\item  Suppose $f$ is distinguished in the variable $y^1_1$, then
  $$\g^f_1= \deg f(0, \ldots, y^1_1 ,\ldots, 0).$$
\item  Suppose $f$ is distinguished in the variable $x_1$, then 
  $$\g^f_0= \deg f(x_1,\ldots, 0).$$
\end{itemize}
\end{prop}

\begin{pf}
{\bf Step 1.}
  We need first to exhibit a suitable basis $c_0^r,\ldots, c_s^r, \ldots, c_r^r$ of
  $A_1 \drk \pn$ (we name equally closed subsets, their associated cycles, and
  their rational classes). These will be the closure in $\drk \pn$ of the one
  dimensional subschemes ${\overset{\circ}{c}}{}_s^r \subseteq U^r$
  defined by the vanishing of all coordinates but $y^1_{1,\stackrel{(s)}{\ldots},{1}}$
  (if $s=0$, this is the coordinate $x_1$). That they form a basis will be proved in Step 2.

  To give an intrinsic description of ${\overset{\circ}{c}}{}_s^r$,
  we define the subspace $\CC^{k-1} \subseteq \CC^k$ in the given affine plane
  $\CC^n \subseteq U^0$ as the space corresponding to coordinates $\seq{x}{2}{k}$,
  and the rank $(k-1)$-subbundle $\tilde X$ of the bundle $X$ on this plane by 
  $\tilde X (P)=P +\CC^{k-1}$, and correspondingly the trivial rank $(k-1)$-subbundle
  $(\tilde{\S}_r)_{U^r} \iso \tilde X_{U^r}$ of $({\S}_r)_{U^r} \iso X_{U^r}$.
  Analogously, if $\CC \subseteq \CC^{n-k}$ is the $y_1$-axis, we define the rank $1$
  subbundle $\tilde Y$ of $Y$ by $\tilde Y (P)= P + \CC$.

  Furthermore, we define inductively $\tilde U^r= \tilde U^{r-1} 
  \x \CC$ inside $U^r =U^{r-1} \x V^r$ by taking $\tilde U^0=U^0$,
  $$
  \tilde U^1= \Hom (X/\tilde X, \tilde Y) \subseteq \Hom (X,Y) =U^1,
  $$
  and assuming that $\tilde U^r$ is already defined, by taking
  $$
  \tilde U^{r+1}= \Hom ((\S_r)_{U^r}/(\tilde{\S}_r)_{U^r}, 
  T\tilde U^r/U^{r-1})= \Hom (U^r \x \CC, \tilde U^r\x\CC)=\tilde U^r\x\CC
  $$ 
  inside $U^{r+1}=\Hom ((\S_r)_{U^r}, TU^r/U^{r-1})$.

  Let $0$ be the origin of $\CC^n$, and let $0_{s-1}=0_{0,s-1}(0)$ be the 
  origin of $U^{s-1}$. Note that $0_{s-1} \in \tilde{U}_{s-1}$.
  Define ${\overset{\circ}{c}}{}^0 =0 + \CC \subseteq U^0$, where $\CC \subseteq
  \CC^{k}$ is the $x_1$-axis. For $s \geq 1$,
  let ${\overset{\circ}{c}}{}^s \subseteq U^s$ be the $1$-dimensional subvariety $\tilde U^s
  \cap b_s^{-1}(0_{s-1})$
  and $c^s$ the closure of ${\overset{\circ}{c}}{}^s$ in $D^s_k \pn$. For
  $0 \leq s\leq r$, the above
  ${\overset{\circ}{c}}{}^r_s \subseteq U^r$ is just $0_{s,r}({\overset{\circ}{c}}{}^s)$,
  thus $b_{r,s}({\overset{\circ}{c}}{}_s^r) = {\overset{\circ}{c}}{}^s$ and
  $b_{r,s}(c_s^r) = c^s$.

{\bf Step 2.}
  Now we want to relate the elements $c_0^r,\ldots, c_s^r, \ldots, c_r^r$ of
  $A_1 \drk \pn$ with the elements $C_0^r,\ldots, C_s^r, \ldots, C_r^r$ of
  $A^1 \drk \pn$.  
  One has ${\overset{\circ}{c}}{}_s^r \iso \CC$ and
  therefore $c_s^r\cap (C_0^r\cup \ldots \cup C_r^r)=c^r_s 
  \backslash U^r=c^r_s \backslash {\overset{\circ}{c}}{}_s^r \iso \PP \backslash \CC$
  consists, set-theoretically, of just one point, which we want now to show not to be in any
  $C^r_t$ for $t< s$ and that this point is in fact the {\em schematic} intersection
  $c^r_s \cap C^r_s$, i.e. that $c^r_s \cdot C^r_s =1$, as rational classes. On the one hand,
  this gives 
  an alternative proof that the $C^r_s$ are linearly independent in lemma~\ref{lem:2}. On the
  other hand, this proves that $c^r_s$ are a basis for $A_1 \drk\pn$, since the
  intersection matrix $(C^r_t \cdot c^r_s)$ is lower-triangular. 

  Let $0 \leq t<s \leq r$. 
  Since $C^r_t=b^{-1}_{r,t}(C^t)$, in order to show that $c^r_s \cap C^r_t=\o$
  it is enough to show that $b_{r,t}(c^r_s) \cap C^t$ is empty.
  By construction, $b_{r,t}({\overset{\circ}{c}}{}_s^r)=0 \in U^t$,
  since $t < s$, thus $b_{r,t}(c_s^r)=0$ is disjoint with $C^t$. 
  This proves that, as rational classes, $c_s^r \cdot C^r_t =0$.

  Again, since $C^r_s=b_{r,s}^{-1}(C^s)$, in order to show that $c^r_s \cdot C^r_s=1$, it is
  enough to show that $b_{r,s}(c^r_s)\cdot C^s=c^s \cdot C^s=1$. For $s=0$ this is true 
  since $c^0 \subseteq \pn$ is a line and $C^0$ is the hyperplane. For $s \geq 1$, we note 
  that since $c^s\subseteq 
  b_s^{-1}(0_{s-1})$,
  it is enough to show for $C^s_{0_{s-1}}=C^s \cap b_s^{-1}(0_{s-1})$ that
  $$
   c^s \cdot C^s_{0_{s-1}}=1
  $$
  in the Grassmannian $b^{-1}_s (0_{s-1})$. Suppose now that $s \geq 2$ (the case $s=1$ is
  similar and is left to the reader). Then $b^{-1}_s (0_{s-1})
  =\Gr_k ((\cF_{s-1})_{0_{s-1}})$, where
  $(\cF_{s-1})_{0_{s-1}}$ is the Semple vector space
  $$
  (\cF_{s-1})_{0_{s-1}}=(TU^{s-1}/U^{s-2})_{0_{s-1}}  \oplus (\S_{s-1})_{0_{s-1}}.
  $$
  In this Grassmannian, $C^s_{0_{s-1}}$ is the Schubert cycle of $k$-spaces meeting 
  $(TU^{s-1}/U^{s-2})_{0_{s-1}}$,
  i.e. the base of its Picard group $A^1$. On the other hand, the $1$-dimensional subvariety
  $$
   {\overset{\circ}{c}}{}^s= 
   \Hom((\S_{s-1})_{0_{s-1}}/(\tilde{\S}_{s-1})_{0_{s-1}}, (T\tilde U^{s-1}/U^{s-2})_{0_{s-1}})
  $$
  of the open subset $\Hom ((\S_{s-1})_{0_{s-1}}, (TU^{s-1}/U^{s-2})_{0_{s-1}})$ of
  $\Gr_k ((\cF_{s-1})_{0_{s-1}})$ is an open subset of the Grassmannian of $k$-subspaces
  of the $(k+1)$-subspace 
  $(T\tilde U^{s-1}/U^{s-2})_{0_{s-1}}  \oplus (\S_{s-1})_{0_{s-1}} \subseteq
  (\cF_{s-1})_{0_{s-1}}$
  which contain the $(k-1)$-subspace $(\tilde{\S}_{s-1})_{0_{s-1}} \subseteq
  (\S_{s-1})_{0_{s-1}} \subseteq (\cF_{s-1})_{0_{s-1}}$. Thus $c^s$ is the Grassmannian
  of such $k$-spaces, i.e. the Schubert cycle base of $A_1$ dual to the base $C^s_{0_{s-1}}$ of
  $A^1$, thus $c^s \cdot C^s_{0_{s-1}}=1$.

{\bf Step 3.} 
  Now we can find the coefficient $\g^f_r$ in the expression of definition~\ref{def:3}.
  By the above, this is
  $$ 
  \g^f_r=[\drk \pn (f)]\cdot c^r_r = \text{length } U^r(f) \cap  {\overset{\circ}{c}}{}^r_r
  = \deg f(0, \ldots, y^1_{1,\stackrel{(r)}{\ldots},{1}} ,\ldots, 0).
  $$
  Indeed, the second equality is due to
  $$
  \left(\drk \pn (f) \backslash U^r(f)\right) \cap  c^r_r
  = \drk \pn (f) \cap \left( C^r_0 \cup \cdots \cup C^r_r \right) \cap  c^r_r=\o
  $$
  since $( C^r_0 \cup \cdots \cup C^r_r) \cap  c^r_r$ is schematically one point, say the
  point $(0,\ldots, \infty, \ldots, 0)$ at infinity of ${\overset{\circ}{c}}{}^r_r=
  \{(0,\ldots, y^1_{1,\stackrel{(r)}{\ldots},1}, \ldots, 0)\} \iso \CC$, and this point
  is not in $\drk\pn (f)$, i.e. 
  $$
   0\neq f(0,\ldots, \infty, \ldots, 0) = \lim_{y^1_{1\ldots1} \ar \infty} 
   f(0,\ldots, y^1_{1,\stackrel{(r)}{\ldots},1}, \ldots, 0),
  $$
  because $f$ is distinguished in the variable $y^1_{1,\stackrel{(r)}{\ldots},1}$.
  This finishes the proof of the first item in the statement of the proposition.

{\bf Step 4.} 
  To get the coefficient $\g^f_1$, we prove first that $c^r_1\cdot C^r_s=0$ for $s>1$. For this
  we consider the (flat)
  family of subvarieties $\{S_{\l}\}_{\l \in \PP^1}$ of $\pn$ such that for $\l \in
  \CC$, $S_{\l}$ is given by the equations $y_1-\l x_1=0$, $y_2=0$, $\ldots$,
  $y_{n-k}=0$, and $S_{\infty}$ has equations $x_1=0$, $y_2=0$, $\ldots$,
  $y_{n-k}=0$. Clearly the points of ${\overset{\circ}{c}}{}^r_1$ are parametrized by
  $D^r_k S_{\l} \cap b_{r,0}^{-1}(0)$, $\l \in \CC$. Therefore $c_1^r \backslash
  {\overset{\circ}{c}}{}^r_1$ is the point $D^r_k S_{\infty} \cap b_{r,0}^{-1}(0)$. But $S_{\infty}$
  is smooth, so $D_k^r S_{\infty}$ is disjoint from $C^r_s=0$ for $s>1$
  (see remark~\ref{rem:smooth}). This proves that 
  $c^r_1\cdot C^r_s=0$ for $s>1$.

  Now $\g^f_1=[\drk \pn (f)]\cdot c^r_1= \deg f(0, \ldots, y^1_1 ,\ldots, 0)$, when $f$ is distinguished
  in the variable $y^1_1$, arguing as in step 3.

{\bf Step 5.} 
  To get the coefficient $\g^f_0$, we need to prove first that $c^r_0\cdot C^r_s=0$ for $s>0$. 
  This time we fix the subvariety $S$ with equations $y_1=0$, $y_2=0$, $\ldots$, $y_{n-k}=0$. 
  Recall ${\overset{\circ}{c}}{}^0 =0 + \CC \subseteq U^0$, where $\CC \subseteq
  \CC^{k}$ is the $x_1$-axis and $c^0$ is its closure in $\pn$.	
  The points
  of $c^r_0$ are parametrized by $D_k^r S \cap b_{r,0}^{-1}(\l)$, $\l \in c^0 \subseteq \pn$. Now
  $D_k^r S$ is disjoint from $C^r_s$, for $s\geq 1$, so $c^r_0\cdot C^r_s=0$ for $s>0$.
  The rest of the argument is as in steps 3 and 4.
\end{pf}

\begin{rem}
\label{rem:c0,c1}
  The explicit computation of the other $\g^f_s$, $2 \leq s <r$, requires the full knowledge of the
  intersection matrix $c_s^r \cdot C_t^r$, which is a more delicate issue (see remark~\ref{rem:delic}). 
  Nonetheless if we are only interested
  in applications of differential equations on smooth varieties, $\g^f_0$ and $\g^f_1$ will suffice 
  (see corollary~\ref{cor:extra}). 
\end{rem}

\begin{rem}
\label{rem:util}
  Suppose that $f(x_i,y_j, y^j_i, y^j_{i_1,i_2}, \cdots, y^j_{i_1,\ldots,i_r})$ is a 
  differential equation on $\pn$ relative to a reference $x_i, y_j$ with
  $d=\deg(f)$. Suppose that 
  $$
   \deg f(0,\ldots, 0, y^1_{1,\stackrel{(r)}{\ldots},1},\ldots, 
   y^{n-k}_{k,\stackrel{(r)}{\ldots},k}) =d.
  $$
  Then $\g_r^f=d$. The proof is similar to steps 2 and 3 in the proof of 
  proposition~\ref{prop:7}. 
  Choose a generic line $\overset{\circ}{l} \subseteq V^r$ and identify it with
  $0_{r-1} \x \overset{\circ}{l} \subseteq U^r=U^{r-1} \x V^r$. 
  Let $l \subseteq \drk \pn$
  be its closure. Then it is easy to see that $l \cdot C_s^r =0$ for $s<r$ and
  $l \cdot C_r^r=1$. Now 
  $$
  \g^f_r =[\drk \pn (f)]\cdot l=\deg f(0,\ldots, 0, y^1_{1,\stackrel{(r)}{\ldots},1},\ldots, 
   y^{n-k}_{k,\stackrel{(r)}{\ldots},k}) =d,
  $$
  as in step 3 in the proof of proposition~\ref{prop:7}. 
\end{rem}

\section{Degree of the divisor of solutions}
\label{sec:3}

Let $S$ be a subvariety of $\pn$ of dimension $k$. We have defined $\drk S
\subseteq \drk \pn$ as the closure of $\drk S_{\reg} \subseteq \drk \pn$.
There is a natural map 
$$
 \drk S \to S,
$$
which is an isomorphism over the 
smooth part. Therefore if $S$ is smooth then $S \iso \drk S$.
In general, we have the following invariants associated to $S$
\begin{defn}
\label{def:cuspidal}
  The cuspidal numbers $\g_S^s$, $s \geq 0$, of $S$ are defined as
\begin{itemize}
  \item $\g^0_S$ is the degree of $S \subseteq \pn$.
  \item $\g^1_S$ is the class of $S$, i.e. the degree of the divisor of $S$ 
       consisting of the points of tangency of tangent $k$-planes to $S$
       that can be drawn from a generic $(n-k-1)$-plane of $\pn$.
  \item For $s \geq 2$, $\g^s_S$ is 
  the degree of the push-forward $b_{s,0} (C^s S) \subseteq \pn$ of 
  $C^s S=D^s_k S \cap C^s \subseteq D^s_k \pn$ under
  $b_{s,0} : D^s_k\pn \ar \pn$.
\end{itemize}
\end{defn}

\begin{rem}
\label{rem:smooth}
  As remarked in~\cite{ASS}, if $S$ is smooth then $\g^s_S=0$ for $s \geq 2$. Indeed, in
  this case, $\drk S$ is disjoint from the cuspidal divisors $C^r_s$, $s \geq 2$.
  Moreover if $S$ has singularities only in codimensions $2$ or more (e.g.\ $S$ normal), 
  the cycle $b_{s,0}
  (C^sS) \subseteq \pn$, $s \geq 2$, has dimension less or equal than $k-2$ and hence 
  $\g^s_S=0$ for $s \geq 2$. 
  If $S$ has no cuspidal singularities (for instance, if the only singularities are normal double
  crossings along a smooth subvariety) then $\g^s_S=0$ for $s\geq 2$ as well. The cuspidal numbers
  measure in some sense how complicated the singularities (on codimension $1$) of $S$ are.
\end{rem}

\begin{lem}
\label{lem:8}
  For any $0 \leq s \leq r$, it is
  $\g^s_S= \drk S \cap C^r_s \cap H^{k-1}$, where $H=C^r_0=b_{r,0}^{-1}(H)$ is the 
  hyperplane in $\drk\pn$.
\end{lem}

\begin{pf}
  Using that 
  $\drk S \cap C^r_s \cap b^{-1}_{r,0}(H^{k-1}) =D^s_k  S \cap C^s \cap 
  b^{-1}_{s,0}(H^{k-1})$, we reduce to the case $r=s$. Now for $r=0$ is obvious, for $r=1$
  follows from the definition, and for $r \geq 2$,
  $\drk S \cap C^r \cap 
  b^{-1}_{r,0}(H^{k-1})= C^r S \cdot b^{-1}_{r,0}(H^{k-1})=
  b_{r,0} (C^r S) \cdot H^{k-1}=\g^r_S$. 
\end{pf}

Given a differential equation $f$ on $\pn$, we define the divisor of solutions of $f$ on $S$
as follows. First, $\drk S(f)$ is the schematic intersection in $\drk \pn$
$$
  \drk S(f) = \drk S \cap \drk \pn (f).
$$
Unless $\drk S(f)=\drk S$, we have that $\drk S(f)$ is a divisor in $\drk S$. 
The divisor of solutions
$S(f) \subset S$ of $f$ on $S$ is the push-forward of $\drk S(f)$ under
$b_{r,0}:\drk \pn \ar \pn$, i.e.\ $S(f)=b_{r,0}(\drk S(f))$.

\begin{thm}
\label{thm:main}
  Let $f$ be a differential equation on $\pn$ with $\g^f_s$, $0 \leq s \leq r$.
  Let $S$ be a subvariety
  of $\pn$ of dimension $k$. Suppose that $S(f)$ is a proper subset of $S$. 
  Then the degree of $S(f)$ is 
  $$
   \g_0^f\g^0_S + \cdots +\g_s^f\g^s_S + \cdots +\g_r^f\g^r_S.
  $$
\end{thm}

\begin{pf}
  The numbers $\g^f_s$ are defined by the condition
  $$
  [\drk \pn (f)]=  \g^f_0 C_0^r +\cdots+ \g^f_s C_s^r+\cdots + \g^f_r C_r^r
  $$
  in $A^1\drk\pn$. Now the degree of $S(f)$ is $S(f) \cdot H^{k-1}$ i.e.
  $$
  \drk S(f) \cdot b_{r,0}^{-1}(H^{k-1})= \drk S \cap \drk \pn (f) \cap H^{k-1}= 
  [\drk S \cap H^{k-1}] \cdot [\drk \pn (f)].
  $$
  Lemma~\ref{lem:8} says that the $s^{th}$-cuspidal degree of $S$ is given as
  $\g^s_S= [\drk S \cap H^{k-1}]\cdot C^r_s$. Therefore we have
  $$
  \deg S(f)= \g_0^f\g^0_S + \cdots +\g_s^f\g^s_S + \cdots +\g_r^f\g^r_S.
  $$
\end{pf}

\begin{rem}
\label{rem:delic}
  Theorem~\ref{thm:main} can be used for computing $\g_s^f$ for a given differential 
  equation $f$, by
  using a standard set of subvarieties whose cuspidal numbers are known (or easily obtainable). 
  This method is used in the examples of section~\ref{sec:4}.
\end{rem}

In the smooth case we get the following

\begin{cor}
\label{cor:extra}
  Let $f$ be a differential equation on $\pn$ and let $S$ be a smooth (or just 
  normal) subvariety
  of $\pn$ of dimension $k$. Suppose that $S(f)$ is a proper subset of $S$. 
  Then  
  $$
   \deg S(f) = \g_0^f\g^0_S + \g_1^f\g^1_S.
  $$
\end{cor}

Let us work out the values of $\g^0_S$ and $\g^1_S$ for a smooth $S$.
In the case $k=1$, $S \subseteq \pn$ is a smooth curve. If $d$ is its degree and
$g$ is genus, then $\g^0_S=d$ and $\g^1_S=2g-2+2d$. Note that when $n=2$, i.e.\
$S\subseteq \PP^2$ is a smooth plane curve, the class of $S$ is $\g^1_S= d(d-1)$, which
can be obtained by using the adjunction formula $2g-2=d(d-3)$.

For $k>1$, let $S\subseteq \pn$ be a smooth (or just normal)
subvariety of degree $d$ and let $g$ be the 
genus of the generic 
section $C=S\cap H^{k-1}$ (which is a smooth curve). Then the degree of $C$ is $d$ and
its class is $\g^1_C=2g-2+2d$. Now it is easy to see that the class of $S$ 
equals that of $C$,
$\g^1_S=\g^1_C$ (for instance take a reference $x_i,y_j$ such that the $H^{k-1}$ has
equations $y_1=0,\ldots, y_{n-k}=0, x_1=0$ and use lemma~\ref{lem:8}).

\begin{cor}
\label{cor:extra2}
  Let $f$ be a differential equation on $\pn$ and let $S \subseteq\pn$ be 
  a smooth (or just normal)
  subvariety of degree $d$ whose generic section $S\cap H^{k-1}$ has genus $g$. 
  Suppose that $S(f)$ is a proper subset of $S$. Then $\deg S(f) = \g_0^f\, d 
  + \g_1^f\, (2g-2+2d)$. Furthermore, if $S$ is a hypersurface, then the formula
  is reduced to $\deg S(f) = \g_0^f\, d + \g_1^f\, d(d-1)$. 
\end{cor}

\section{Examples}
\label{sec:4}

Theorem~\ref{thm:main} can be used in two directions, either to compute the degree of the 
divisor of solutions of a differential equations on a subvariety of $\pn$
or to extract information
about a differential equation $f$ and a subvariety $S\subseteq\pn$ once we have computed
the degree of the divisor of solutions of $f$ on $S$.

\noindent {\bf Computation of $\text{deg } S(f)$.}
Given a differential equation $f$ and a $k$-dimensional subvariety $S \subseteq \pn$,
to find the divisor $S(f)$ we proceed as follows. Take a reference $x_i,y_j$. When $S$ is
smooth we find explicitly $U^r(f) \cap \drk S$ by considering the ideal generated by 
the equations defining $S$, their formal derivatives up to order $r$ together with the
equation $f=0$. If $S$ is non-smooth we have to be a little bit more careful, as
$\drk S$ is defined as the closure of $\drk S_{\reg}$, and hence it may be smaller 
than the set defined by the above ideal (see example~\ref{ex:b}).

Using different references $x_i,y_j$, we find with this method all the points in $\drk S(f)$ 
not lying in 
\begin{equation}
\label{eqn:sols4}
  (C^r_2 \cup \ldots \cup C^r_r) \cap \drk \pn (f)\cap \drk S,
\end{equation}
i.e. $(\drk S(f))_{\text{nc}} =\drk S(f) - \bigcup\limits_{s=2}^{r} C^r_s$. When $S$ is
smooth \eqref{eqn:sols4} is empty and hence $\drk S(f)=(\drk S(f))_{\text{nc}}$.
In general,
$\deg S(f)$ equals the degree of the closure of $(\drk S(f))_{\text{nc}}$
unless~\eqref{eqn:sols4} has $(k-1)$-dimensional components (we expect it to be
$(k-2)$-dimensional).

Note that to compute $\deg S(f)$ we have to take a linear section $H^{k-1}$ of $S(f)$. Therefore
we look for the number of solutions of the differential equation in a general $H^{k-1}$ 
section of $S$, which is a curve. 

\begin{ex}
\label{ex:a}
Let $n=2$ and $k=1$, i.e.\ the case of curves in $\PP^2$.  Consider the smooth cubic
$S$ given by $x^3+y^2=1$. To compute $\deg S(y')$, where $y'={\bd y \over \bd x}$,
we compute the number of points in
$$
  \left\{ \begin{array}{l} x^3 +y^2 =1 \\ 3x^2+2yy' =0 \\ y'=0 \end{array}\right.
$$
which is $4$. Now $S$ passes through the point $P=(0,\infty)$ at infinity $C^0$ and
is tangent to $C^0$ there. Therefore this point counts with multiplicity $2$ (otherwise
compute in a different reference). So $\deg S(y')=4+2=6$.

Now let us compute $\deg S(y'')$, where $y''={\bd^2 y \over \bd x^2}$.
The number of points in
$$
  \left\{ \begin{array}{l} x^3 +y^2 =1 \\ 3x^2+2yy' =0 \\ 6x+2(y')^2+2yy''=0 \\ y''=0 
  \end{array}\right.
$$
is $8$. The point $P$ at infinity counts once, since it is a simple flex of $S$. Hence
$\deg S(y'')=8+1=9$.
\end{ex}

\begin{ex}
\label{ex:b}
Now consider the singular cubic $T$ given by $x^3+y^2=0$. Now the closure of 
$D_1^1 T_{\reg}$ is the irreducible component of
$$
  \left\{ \begin{array}{l} x^3 +y^2 =0 \\ 3x^2+2yy' =0 \end{array}\right.
$$
given by the equations $x={2\over 3}(y')^2$, $y=-{2\over 3}(y')^3$. Hence $U^0 \cap T(y')$
consists of $1$ point. Also $T$ passes through $P=(0,\infty)$ 
at infinity with multiplicity $2$, so $\deg T(y')=1+2=3$.

For computing $\deg T(y'')$ we study
$$
  \left\{ \begin{array}{l} 3x-2(y')^2=0 \\ 3y+2(y')^3=0 \\ 3-4y'y'' =0 \\ 3y'+6(y')^2y''=0 
  \\ y''=0 \end{array}\right.
$$
which is empty. Counting the point at infinity, $\deg T(y'')=1$.
\end{ex}

\noindent {\bf Computation of $\g_s^f$.}
Apart from proposition~\ref{prop:7} (and in some cases the extension given in
remark~\ref{rem:util}), we can compute $\g_s^f$ by coupling $f$ with some simple 
examples of varieties $S \subseteq \pn$.
Many differential equations $f$ have a geometrical meaning. For instance for curves in
the plane, $y'=0$ gives the points of a curve with tangent parallel
to the $x$-axis, and $y''=0$ gives the
flexes of a curve. Let us work out these two examples.

For $f=y'$, we have $\g^{y'}_1=1$. Let $C \subseteq \PP^2$ be a smooth conic, which has degree
$\g^0_C=2$ and class $\g^1_C=2$. Hence $2=\deg C(y')= 2\g^{y'}_0 +2$, so $\g^{y'}_0=0$. 

For $f=y''$, proposition~\ref{prop:7} says that $\g^{y''}_2=1$. A conic has no flexes, so
$0=\deg  C(y'')= 2\g^{y''}_0+2 \g^{y''}_1$. The smooth cubic $S$ of example~\ref{ex:a} gives
first $\g^S_1=\deg S(y')=6$ and then 
$9=\deg  S(y'')= 3\g^{y''}_0+6 \g^{y''}_1$. Therefore $\g^{y''}_0=-3$ and 
$\g^{y''}_1=3$. 

\noindent {\bf Computation of $\g^s_S$.}
For a smooth (or normal)
variety $S \subseteq \pn$, we have that $\g^0_S$ is its degree and $\g^1_S$
its class. This can be computed by taking a general section $C=H^{k-1}\cap S$, which
is a curve of degree $\g^0_C=\g^0_S$ and class $\g^1_C=\g^1_S$. So $\g^1_S=2g-2+2d$, where
$g$ is the genus of $C$ and $d$ its degree. Also $\g^s_S=0$ for $s \geq 2$.

In the non-smooth case, things are a little bit more complicated. For instance, for
the smooth cubic $S$ of example~\ref{ex:a}, $\g^0_S=3$, $\g^1_S=6$ and $\g^2_S=0$.
Instead for the non-smooth cubic $T$ of example~\ref{ex:b},
$\g^0_T=3$, $\g^1_T=\deg T(y')=3$ and $\g^2_T=\deg T(y'') -\g_0^{y''} 3 -\g_1^{y''} 3=1$.

\noindent {\bf Degree of the divisor of parabolic points.}
As an application, we shall determine the degree of the divisor of parabolic points
of a subvariety $S\subseteq \pn$. 
Suppose first that $k=n-1$, i.e.\ $S$ is a hypersurface of $\pn$.
Parabolic points are those points of $S$ with higher contact with the tangent space than
expected. In terms of a reference $\seq{x}{1}{n-1},y$, they are the solutions to the 
differential equation
$$
   f= \det ({\bd^2 y \over \bd x_i \bd x_j}).
$$
To compute the invariants of $f$ we work as follows. By remark~\ref{rem:util}, $\g_2^f=
n-1$. Now for a smooth quadric $C \subseteq \pn$ there are no parabolic points, so 
$0=\deg C(f)=2\g^f_0 +2\g^f_1$, and
$\g^f_1=-\g^f_0$. Now let $S\subseteq \pn$ be the smooth cubic given by
$x_1^3+\cdots +x_{n-1}^3 +y^2 +1=0$. It is easy to compute
$$
 \left\{  \begin{array}{ll} {\bd^2 y\over \bd x_i^2}= -{3x_i \over y}-{9 x_i^4 \over 4y^3}, 
    \qquad & 1 \leq i \leq n-1 \\
  {\bd^2 y\over \bd x_i\bd x_j}= {-9 x_i^2x_j^2 \over 4y^3}, & 1 \leq i,j \leq n-1, \, i\neq j
   \end{array} \right.
$$
so that $\det ({\bd^2 y \over \bd x_i \bd x_j})= (-3)^{n-1}x_1 \cdots x_{n-1}
 (y^2-3) / 4y^{n+1}$ and hence $3(n+1)=\deg S(f) =3\g^f_0 + 6\g^f_1$. This yields
$\g^f_1=-\g^f_0=n+1$. Our conclusion is that for a hypersurface $S \subseteq \pn$, the
degree of the divisor of parabolic points is
$$
  \deg S(f)= -(n+1) \g^0_S + (n+1) \g^1_S +(n-1) \g^2_S.
$$
If $S$ is a smooth hypersurface of degree $d$, it reduces to
$$
  \deg S(f)= -(n+1) d + (n+1) d(d-1)= (n+1) d (d-2).
$$
This agrees with the following alternative argument (only applies 
to the smooth case):
if $S$ is given by the equation $F(\seq{x}{0}{n})=0$, 
then the parabolic points are the intersection of
$F=0$ and $\det({\bd^2 F \over \bd X_i\bd X_j})=0$, and so form
a divisor of degree $(n+1)d(d-2)$.

For the general case $0 < k<n$, we consider a reference $x_i, y_j$ 
and fix $1 \leq r \leq n-k$. Then the parabolic points in the $y_r$-direction
are the solutions to the differential equation
$$
   f_r= \det ({\bd^2 y_r \over \bd x_i \bd x_j}).
$$
Working as above we have the formula 
$$
  \deg S(f_r)= -(n+1) \g^0_S + (n+1) \g^1_S +(n-1) \g^2_S.
$$

\section{The real case}
\label{sec:5}

In this section, our purpose is to extend theorem~\ref{thm:main} to the case of real projective
varieties inside the real projective space $\pnr$. We are 
not going to develop the general theory of data schemes for smooth
real algebraic varieties, but only to outline the construction of $\drk\pnr$.

First of all fix a trivialisation on $\pn$ as in section~\ref{sec:2}, so $\pn=\PP(\CC^{n+1})$.
Again $1 \leq k \leq n-1$.
There is an anti-holomorphic involution $\s:\pn \ar \pn$ coming from conjugation on $\CC^{n+1}$,
whose fixed point set is $\pnr$. Inductively, 
$\s$ induces anti-holomorphic involutions $\s_r: \drk \pn \ar \drk \pn$, $r\geq 1$. We define
$\drk \pnr$ to be the fixed point set of $\s_r$. The maps $b_{r,s}$ 
of~\eqref{eqn:sols2} restrict to $\drk \pnr$ and so there are well-defined
maps $b_{r,s}: D^r_k \pnr \ar D^s_k \pnr$.
It is easy to see that we can construct
$\drk \pnr$ as follows. We take $D^0_k \pnr=\pnr$, $D^1_k \pnr=\Gr_k T\pnr$, the real 
Grassmannian of $k$-planes in the tangent bundle $T\pnr$, and for $r \geq 2$, $\drk
\pnr=\Gr_k \cF^{\RR}_{r-1}$, where the real (Semple) bundle $\cF^{\RR}_{r-1}$ is defined
by a diagram as~\eqref{eqn:sols1}, where the relative 
tangent bundles and the universal sequence of
the right hand side are understood to be real.
These $\drk \pnr$ are also smooth compact differentiable manifolds.

The involution $\s_r$ takes $C^r$ to itself and the fixed point set will be called
$C^r_{\RR}$. This is a smooth $1$-codimensional real algebraic subvariety of $\drk \pnr$.
Clearly $C^0_{\RR} \subseteq \pnr$ is the hyperplane. The other $C^r_{\RR}$ can be also
defined as the cuspidal locus of $\drk \pnr$, namely the Schubert special cycle of real
$k$-planes of $\cF^{\RR}_{r-1}$ meeting $TD^{r-1}_k \pnr/D^{r-2}_k \pnr$.
Again $C^r_{\RR,s}$ are defined either as the fixed point set of $\s_r$ on $C^r_s$ or as
$b_{r,s}^{-1}(C^s_{\RR})$. 

We can parallel the discussion in section~\ref{sec:2} to see that 
$$
  U^r_{\RR}=\drk \pnr \backslash (C^r_{\RR,0} \cup \ldots \cup C^r_{\RR,s}
  \cup \ldots\cup C^r_{\RR,r})
$$
are cartesian powers of $\RR$. Indeed, $U^0_{\RR}= \RR^{k} \x \RR^{n-k}$ has real 
coordinates (the restriction to $U^0_{\RR} \subseteq U^0$ of)
$\seq{x}{1}{k},\seq{y}{1}{n-k}$.
In general 
$$
  U^r_{\RR}\iso U^{r-1}_{\RR} \x V^r_{\RR},
$$
where $V^r_{\RR}=\Hom ((\RR^k)^{\ox r}, \RR^{n-k})$, for any $r \geq 1$.
Alternatively, $\s_r$ restricts to 
$U^r$ and the fixed point locus is $U^r_{\RR}$.
The coordinates for $U^{r}_{\RR}$ will be
\begin{equation}
  x_1, \ldots, x_k,y_1, \ldots, y_{n-k} \qquad \text{and} \qquad
   y^j_{i_1,\ldots,i_s}={\bd^s y_j \over \bd x_{i_s} \cdots \bd x_{i_1}},
\end{equation}
for any $1 \leq s \leq r$, $1 \leq j \leq n-k$, 
and $1 \leq \seq{i}{1}{s} \leq k$.

\begin{lem}
\label{lem:r0}
  $C^r_{0,\RR}, \ldots,C^r_{s,\RR}, \ldots, C^r_{r,\RR}$ form a basis for $H^1 (\drk \pnr
  ;\ZZ/2\ZZ)$, where we name equally the algebraic subvarieties and the cohomology classes they 
  represent through Poincar\'e duality.
\end{lem}

\begin{pf}
  Given the description of the real Grassmannian as a homogeneous space
  $$
     \Gr_k \RR^n= O(n)/O(k)\x O(n-k),
  $$
  it is easy to prove that the fundamental group $\pi_1(\Gr_k \RR^n)=\ZZ/2\ZZ$ for $n>2$ and
  $\pi_1(\Gr_1 \RR^2)=\pi_1(\RR\PP^1)=\ZZ$. So $H^1(\Gr_k \RR^n;\ZZ/2\ZZ)=\ZZ/2\ZZ$ and it
  is generated by the cuspidal subvariety. The Serre spectral sequence of the
  fibration $$\Gr_k \RR^{\bullet} \ar \drk\pnr=\Gr_k\cF_{r-1}^{\RR} \ar D^{r-1}_k \pnr$$ 
  implies that $H^1(\drk \pnr;\ZZ/2\ZZ)$ is generated by 
  $C^r_{0,\RR}, \ldots, C^r_{r,\RR}$. Now the $1$-cycles $c^r_{s,\RR}$ are defined as
  in the proof of proposition~\ref{prop:7} (alternatively as the fixed locus of $\s_r$ on $c^r_s$) and
  the homology classes they represent in $H_1(\drk \pnr;\ZZ/2\ZZ)$ satisfy
  $c^r_{t,\RR} \cdot C^r_{s,\RR}=0$, for $t<s$ and $c^r_{s,\RR} \cdot C^r_{s,\RR}=1 \pmod 2$.
  This shows that $C^r_{s,\RR}$ are linearly independent.
\end{pf}

Recall that the cohomology ring of $\pnr$ (with $\ZZ/2\ZZ$-coefficients) 
is $(\ZZ/2\ZZ)[H]/H^{n+1}$, where
$H$ stands for the hyperplane class. Again $H=C^r_0=b_{r,0}^{-1}(H)$ is the 
hyperplane in $\drk\pn$.
Now we are in the position of defining differential equations on $\pnr$.
\begin{defn}
\label{def:r1}
  A differential equation on $\pnr$ relative to a reference $x_i,y_j$ is a
  non-zero algebraic equation with real coefficients
  $$
    f(x_i,y_j, y^j_i, y^j_{i_1,i_2}, \cdots, y^j_{i_1,\ldots,i_r})=0,
  $$ 
  on $U^r_{\RR}$, for some $r \geq 1$. 
\end{defn}

  This time the zero locus of such $f$ is not necessarily a hypersurface of $\drk \pnr$ (it
  might even be empty). But we may still define the numbers $\g^f_s \in \ZZ/2\ZZ$ by looking at 
  the component in $H^1(\drk\pn;\ZZ/2\ZZ)$ defined by the closure of the zero locus of $f$
  in $U^r_{\RR} \subseteq \drk\pnr$. So 
 $$
   [\drk\pnr(f)]_1= \g_0^f C^r_{0,\RR} +\cdots +\g_r^f C^r_{r,\RR}. 
 $$
  We leave the following analogue of proposition~\ref{prop:7} to the reader.
\begin{prop}
  Let $f$ be a differential equation on $\pnr$ relative to a reference $x_i,y_j$.
\begin{itemize}
\item  If $f$ is distinguished in the variable $y^1_{1,\stackrel{(r)}{\ldots},{1}}$, then
  $\g^f_r= \deg f(0, \ldots, y^1_{1,\stackrel{(r)}{\ldots},{1}} ,\ldots, 0)$.
\item  If $f$ is distinguished in the variable $y^1_1$, then
  $\g^f_1= \deg f(0, \ldots, y^1_1 ,\ldots, 0)$.
\item  If $f$ is distinguished in the variable $x_1$, then 
  $\g^f_0= \deg f(x_1,\ldots, 0)$.
\end{itemize}
\end{prop}

Let $S \subseteq \pnr$ be a real algebraic $k$-dimensional subvariety (this is the zero
locus of polynomial equations with real coefficients such that it has a dense open subset that
is a smooth differentiable manifold of dimension $k$). Considering the same equations 
in $\pn$, we get a complex subvariety $S_{\CC} \subseteq \pn$ on which $\s$ acts with 
fixed point set $S$. We define $\drk S \subseteq \drk \pnr$ as the fixed point set of
$\s_r$ on $\drk S_{\CC} \subseteq \drk\pn$. There is a natural map
$\drk S \to S$, which is an isomorphism over the smooth part.
When $S$ is smooth then $S \iso \drk S$ and $\drk S$ is disjoint from the cuspidal 
divisors $C^r_{s,\RR}$, $s \geq 2$. By analogy with lemma~\ref{lem:8}, we define the
cuspidal numbers of $S$ as follows

\begin{defn}
\label{def:8}
  For any $0 \leq s \leq r$, we define the $s^{th}$-cuspidal number of $S$ as
  $\g^s_S= \drk S \cap C^r_s \cap H^{k-1} \in \ZZ/2\ZZ$ (computed in the cohomology ring
  of $\drk \pnr$).
\end{defn}

Note also that when $S \subseteq \pnr$ is a smooth differentiable $k$-dimensional
manifold, $\drk S \subseteq \pnr$ is also defined by mimicking the algebraic construction 
and $\drk S \to S$ is a diffeomorphism. In this case $\g^0_S$ and $\g^1_S$ are defined
as above.

Given a differential equation $f$ on $\pn$, we define $S(f)$ as the homology class
in $H_{k-1}(\pnr;\ZZ/2\ZZ)$ given as the push-forward of
$$
  \drk S(f) = \drk S \cap \drk \pnr (f) \in H_{k-1}(\drk\pnr;\ZZ/2\ZZ).
$$
Its degree is 
$$
  \deg S(f)=S(f) \cap H^{k-1} \in \ZZ/2\ZZ.
$$ 
The following result is proved much in the same way as theorem~\ref{thm:main}.

\begin{thm}
\label{thm:main2}
  Let $f$ be a differential equation on $\pnr$ with $\g^f_s \in \ZZ/2\ZZ$, $0 \leq s \leq r$.
  Let $S$ be a real algebraic subvariety of $\pnr$ of dimension $k$. Then 
  $$
   \deg S(f)= \g_0^f\g^0_S + \cdots +\g_s^f\g^s_S + \cdots +\g_r^f\g^r_S \in \ZZ/2\ZZ.
  $$
  If $S$ is smooth or $S$ is just a differentiable manifold of dimension 
  $k$, we have
  $$ 
   \deg S(f)= \g_0^f\g^0_S +\g_1^f\g^1_S \in \ZZ/2\ZZ.
  $$ 
\end{thm}

\noindent {\bf Umbilical points.}
  We shall compute the parity of the degree of the subset of umbilical
  points for surfaces in $\PP_{\RR}^3$. For a reference $x_1,x_2,y$, these are
  the points where the Hessian $\left({\bd^2 y \over \bd x_i \bd x_j}\right)$ is diagonal.
  Therefore we look for the solutions to
  $$
   f=\left({\bd^2 y \over \bd x_1^2}- {\bd^2 y \over \bd x_2^2}\right)^2 +4 
   \left({\bd^2 y \over \bd x_1 \bd x_2}\right)^2.
  $$
  Clearly $\g_2^f=2=0$ (we work over $\ZZ/2\ZZ$). For the cubic $S$ given by
  $x_1^3+y^2+1=0$ with $\g^0_S=1$ and $\g^1_S=0$ (see example~\ref{ex:a}),
  the equation $f$ reduces to $({\bd^2 y \over \bd x_1^2})^2$, so 
  $\g_0^f = \deg S(({\bd^2 y \over \bd x_1^2})^2)=2 \deg ({\bd^2 y \over \bd x_1^2})=0$.
  For the cubic $T$ given by
  $x_1^3+y^2=0$, with $\g^0_S=1$ and $\g^1_S=1$ (see example~\ref{ex:b}),
  the equation $f$ reduces to $({\bd^2 y \over \bd x_1^2})^2$, so 
  $\g_0^f +\g_1^f = \deg S(({\bd^2 y \over \bd x_1^2})^2)=0$.
  Therefore $\g_0^f =\g_1^f=\g^f_2=0$ and the number of umbilical points of 
  any surface $S
  \subseteq  \PP_{\RR}^3$ (smooth or not) is always even.

\vspace{5mm}
Departamento de \'Algebra, Geometr\'{\i}a y Topolog\'{\i}a \\
\indent Facultad de Ciencias \\
\indent Universidad de M\'alaga \\
\indent Campus de Teatinos, s/n \\
\indent 29071 M\'alaga \\
\indent Spain

\vspace{1mm}
{\em E-mail:\/} vmunoz@@agt.cie.uma.es

\vspace{5mm}
Departamento de \'Algebra \\
\indent Facultad de Ciencias Matem\'aticas \\
\indent Universidad Complutense de Madrid \\
\indent 28040 Madrid \\
\indent Spain

\vspace{1mm}
{\em E-mail:\/} sols@@mat.ucm.es

\end{document}